\documentclass{amsproc}
\usepackage{amssymb}

\usepackage{amsmath}
\usepackage{graphicx}

\theoremstyle{plain}

\newtheorem{conjecture}{Conjecture}

\newtheorem{definition}{Definition}

\newtheorem{theorem}{Theorem}
\numberwithin{equation}{section}

\begin{document}
\title[Biseparating maps]{Biseparating maps between operator algebras}
\author{Jes\'{u}s Araujo}
\address{Departamento de Matem\'{a}ticas, Estad\'{\i}stica y Computaci\'{o}n\\
Universidad de Cantabria\\
Facultad de Ciencias\\
Avda. de los Castros, s. n.\\
E-39071 Santander, Spain}
\email{araujoj@unican.es\qquad\ http://www.matesco.unican.es/\symbol{126}araujo/ }
\author{Krzysztof Jarosz}
\address{Department of Mathematics and Statistics\\
Southern Illinois University Edwardsville\\
IL 62026, USA}
\email{\ kjarosz@siue.edu\qquad http://www.siue.edu/\symbol{126}kjarosz/}
\thanks{Research of J. Araujo partially supported by the Spanish Direcci\'{o}n
General de Investigaci\'{o}n Cient\'{i}fica y T\'{e}cnica (DGICYT, PB98-1102)}

\begin{abstract}
We prove that a biseparating map between spaces $B\left( E\right) $, and
some other Banach algebras, is automatically continuous and an algebra
isomorphism.
\end{abstract}

\maketitle

\section{Introduction}

Linear maps between Banach algebras, Banach lattices, or Banach spaces
preserving certain properties have been of a considerable interest for many
years. The most classical question concerns isometries, that is, maps that
preserve the norm. More recently, maps that preserve spectrum, spectral
radius, commutativity, normal elements, self-adjoint elements, nilpotents,
idempotents, linear rank, disjointness, or other properties have been
intensely investigated, see for example \cite{BrsSem1993, BrsSem1997, GLS,
JafSou, MarMoy, Molnar,Sem1993, Sou1996, Wolff} and the references given
there. Here we study the biseparating maps, that is maps preserving
disjointness in both directions.

\begin{definition}
A linear map $T$ between algebras $\mathcal{A}$,$\mathcal{B}$ is called
separating if 
\begin{equation*}
ab=0\Rightarrow T\left( a\right) T\left( b\right) =0\text{, \qquad for all }%
a,b\in\mathcal{A}\text{;}
\end{equation*}
it is called biseparating if $T^{-1}:\mathcal{B}\rightarrow\mathcal{A}$ \
exists and is also separating.
\end{definition}

Any algebra isomorphism; that is, a map that preserves both linear and
multiplicative structures, is clearly separating, also any algebra
isomorphism followed by a multiplication by a fixed element from the
commutant of $\mathcal{B}$ is separating - we will call such a map a \emph{%
standard separating map}. However in general a separating map may be very
far from being multiplicative. For example if $\mathcal{A}$ is the disc
algebra then the product of two elements from $\mathcal{A}$ is equal to zero
only if one of them is already zero, consequently any linear map on $%
\mathcal{A}$ is separating. On the other hand any biseparating map $%
T:C\left( X\right) \rightarrow C\left( Y\right) $, where $C\left( X\right) $
is the space of all continuous functions defined on a completely regular set 
$X$, is of the form 
\begin{equation}
Tf=\tau \cdot f\circ \varphi \text{, \qquad for all }f\in C\left( X\right) ,
\label{can}
\end{equation}
where $\tau $ is a nonvanishing continuous function and $\varphi $ is a
homeomorphism from $rY$, the realcompactification of $Y$, onto $rX$ \cite
{ABN}. It is interesting to notice that for compact sets $X,Y$ a separating
invertible map $T:C\left( X\right) \rightarrow C\left( Y\right) $ is
automatically biseparating and consequently of the form (\ref{can}) \cite
{Jar1990}. However, whether this is also true for noncompact sets $X,Y,$ is
an open problem (see \cite{AraJar} for a partial solution). A reader
interested in separating and biseparating maps in more general setting may
want to check a recent monograph \cite{AbrKit2000}.

In this note we show that any biseparating map between the algebras $B\left(
E\right) $ of all continuous linear maps on a Banach space $E,$ as well as
between certain subalgebras of $B\left( E\right) $ and tensor products of
such algebras is standard. The results apply both in the real and in the
complex cases.

\section{Results}

For a Banach space $E$ we will call a subalgebra $\mathcal{A}$ of $B\left(
E\right) $ \emph{standard} if it contains all finite-dimensional operators
and the identity operator $Id$. Some authors also assume that a standard
subalgebra is closed in the norm topology, here we do not make this
assumption. Notice that a standard subalgebra must contain all continuous
projections onto closed finite-codimensional subspaces of $E$. Indeed if $P$
is such a projection then $P=Id-\left( Id-P\right) ,$ where $Id-P$ is finite
dimensional.

\begin{theorem}
\label{th1}Assume $E_{1},E_{2}$ are Banach spaces and $\mathcal{A}_{1},%
\mathcal{A}_{2}$ are standard subalgebras of $B\left( E_{1}\right) $ and of $%
B\left( E_{2}\right) $, respectively. If $T:\mathcal{A}_{1}\rightarrow 
\mathcal{A}_{2}$ is biseparating then it is continuous and of the form 
\begin{equation*}
T\left( A\right) =\alpha S\circ A\circ S^{-1}\text{,\quad\ for }A\in 
\mathcal{A}_{1},
\end{equation*}
where $S$ is a continuous linear isomorphism from $E_{1}$ onto $E_{2}$ and $%
\alpha $ a nonzero scalar.
\end{theorem}

\begin{proof}
We first need to introduce some notation. For $i=1,2$ we put 
\begin{align*}
B^{\bot }& =\{A\in \mathcal{A}_{i}\backslash \{0\}:A\circ B=0\}\text{, \quad
for }B\in \mathcal{A}_{i}\text{,} \\
\mathcal{M}_{i}& =\{B^{\bot }:B\in \mathcal{A}_{i}\}, \\
\mathcal{M}_{i}^{1}& =\mathcal{M}_{i}\backslash \{B^{\bot }\in \mathcal{M}%
_{i}:\exists C\in \mathcal{A}_{i},\text{ }\emptyset \neq C^{\bot
}\varsubsetneq B^{\bot }\}.
\end{align*}
For $e\in E_{i}$ and $e^{\ast }\in E_{i}^{\ast }$ we denote by $e\otimes
e^{\ast }$ the one dimensional operator on $E_{i}$ defined by 
\begin{equation*}
\left( e\otimes e^{\ast }\right) \left( e^{\prime }\right) =e^{\ast }\left(
e^{\prime }\right) e,\text{\quad\ for }e^{\prime }\in E_{i}.
\end{equation*}
Since the definitions of $\mathcal{M}_{1}$ and $\mathcal{M}_{1}^{1}$ involve
only the structures that are preserved by $T$, $T$ maps these sets on $%
\mathcal{M}_{2}$ and $\mathcal{M}_{2}^{1}$, respectively.

Notice that $B^{\bot }$ is trivial if and only if $\mathrm{range}\left(
B\right) $ is dense, and $B^{\bot }$ consists of operators of dimension one
if and only if the closure $\overline{\mathrm{range}B}$ of the $\mathrm{range%
}B$ is one-codimensional. Since $\mathcal{A}_{i}$ contains all projections
onto closed one-codimensional subspaces, we have 
\begin{equation*}
\mathcal{M}_{i}^{1}=\{\{A\in \mathcal{A}_{i}:\ker A=E_{0}\}:E_{0}=\overline{%
E_{0}}\subset E_{i}\text{, }\dim E_{i}/E_{0}=1\}.
\end{equation*}
Hence $\bigcup \mathcal{M}_{i}^{1}$ is simply equal to the set of all the
one dimensional operators, so $T$ maps a one dimensional operator onto a one
dimensional operator.

Fix a linear continuous functional $e^{\ast }$ on $E_{1}$. For any $e\in
E_{1}$ we have 
\begin{equation*}
T\left( e\otimes e^{\ast }\right) =f\otimes f_{e^{\ast }}^{\ast }
\end{equation*}
for some $f\in E_{2}$ and $f_{e^{\ast }}^{\ast }\in E_{2}^{\ast }$. If we
change the point $e$ but keep the same functional $e^{\ast }$, the operator $%
e_{1}\otimes e^{\ast }$ will still belong to the same element of $\mathcal{M}%
_{1}^{1}$. Hence $T\left( e_{1}\otimes e^{\ast }\right) =f_{1}\otimes
f_{1}^{\ast }$ belongs to the same element of $\mathcal{M}_{2}^{1}$; this
means that $\ker f_{1}^{\ast }=\ker f_{e^{\ast }}^{\ast }$ and the
functionals $f_{e^{\ast }}^{\ast },f_{1}^{\ast }$ are proportional. So there
is a map $S_{e^{\ast }}:E_{1}\rightarrow E_{2}$ such that 
\begin{equation*}
T\left( e\otimes e^{\ast }\right) =S_{e^{\ast }}\left( e\right) \otimes
f_{e^{\ast }}^{\ast }\text{, for all }e\in E_{1}.
\end{equation*}
Since $T$ is linear, so must be $S_{e^{\ast }}$; since $\ker T$ is trivial,
the same must be true about $\ker S_{e^{\ast }}$. Finally, since $T$ is
surjective, it maps an element of $\mathcal{M}_{1}^{1}$ onto an entire
element of $\mathcal{M}_{2}^{1}$ so $S_{e^{\ast }}$ is surjective.

Assume $e_{1}^{\ast }$ is another continuous linear functional on $E_{1}$,
not proportional to $e^{\ast },$ and let $S_{e_{1}^{\ast }},f_{e_{1}^{\ast
}}^{\ast }$ be such that 
\begin{equation*}
T\left( e\otimes e_{1}^{\ast }\right) =S_{e_{1}^{\ast }}\left( e\right)
\otimes f_{e_{1}^{\ast }}^{\ast }\text{,\quad\ for all }e\in E_{1}.
\end{equation*}
Since $\ker e^{\ast }\neq \ker e_{1}^{\ast },$ then $e\otimes e^{\ast }$ and 
$e\otimes e_{1}^{\ast }$ belong to distinct elements of $\mathcal{M}_{1}$,
for any $e\in E_{1}\backslash \{0\}$, so $\ker f_{e^{\ast }}^{\ast }\neq
\ker f_{e_{1}^{\ast }}^{\ast }.$ Suppose the linear maps $S_{e^{\ast }}$ and 
$S_{e_{1}^{\ast }}$ are not proportional and let $e_{0}$ be such that $%
S_{e^{\ast }}\left( e_{0}\right) $ and $S_{e_{1}^{\ast }}\left( e_{0}\right) 
$ are linearly independent. Then the operator 
\begin{equation*}
e_{0}\otimes e^{\ast }+e_{0}\otimes e_{1}^{\ast }
\end{equation*}
is one dimensional, while 
\begin{equation*}
T\left( e_{0}\otimes e^{\ast }+e_{0}\otimes e_{1}^{\ast }\right) =S_{e^{\ast
}}\left( e_{0}\right) \otimes f_{e^{\ast }}^{\ast }+S_{e_{1}^{\ast }}\left(
e_{0}\right) \otimes f_{e_{1}^{\ast }}^{\ast },
\end{equation*}
is two dimensional. The contradiction shows that there is a linear bijection 
$S:E_{1}\rightarrow E_{2}$ and a map $\Psi :E_{1}^{\ast }\rightarrow
E_{2}^{\ast }$ such that 
\begin{equation}
T\left( e\otimes e^{\ast }\right) =S\left( e\right) \otimes \Psi \left(
e^{\ast }\right) \text{, \quad for all }e\in E_{1},e^{\ast }\in E_{1}^{\ast
}.  \label{2}
\end{equation}
As\ in the case of $S_{e^{\ast }}$, since $T$ is a linear bijection so must
be $\Psi $.

We now show that $S$ is continuous. For any $e_{1},e_{2}\in E_{1}$ and $%
e_{1}^{\ast },e_{2}^{\ast }\in E_{1}^{\ast }$ we have 
\begin{equation*}
\left( e_{1}\otimes e_{1}^{\ast }\right) \circ \left( e_{2}\otimes
e_{2}^{\ast }\right) =e_{1}^{\ast }\left( e_{2}\right) \left( e_{1}\otimes
e_{2}^{\ast }\right) ,
\end{equation*}
and by (\ref{2}) 
\begin{align*}
T\left( e_{1}\otimes e_{1}^{\ast }\right) \circ T\left( e_{2}\otimes
e_{2}^{\ast }\right) & =\left( S\left( e_{1}\right) \otimes \Psi \left(
e_{1}^{\ast }\right) \right) \circ \left( S\left( e_{2}\right) \otimes \Psi
\left( e_{2}^{\ast }\right) \right) \\
& =\left( \Psi \left( e_{1}^{\ast }\right) \right) \left( S\left(
e_{2}\right) \right) \left( S\left( e_{1}\right) \otimes \Psi \left(
e_{2}^{\ast }\right) \right) .
\end{align*}
So 
\begin{equation*}
e_{1}^{\ast }\left( e_{2}\right) =0\text{ iff }\left( \Psi \left(
e_{1}^{\ast }\right) \right) \left( S\left( e_{2}\right) \right) =0,
\end{equation*}
hence 
\begin{equation*}
\ker \Psi \left( e_{1}^{\ast }\right) =S\left( \ker e_{1}^{\ast }\right) 
\text{,\quad\ for all }e_{1}^{\ast }\in E_{1}^{\ast }
\end{equation*}
or 
\begin{equation}
\ker \left( \Psi \left( e_{1}^{\ast }\right) \circ S\right) =S^{-1}\left(
\ker \Psi \left( e_{1}^{\ast }\right) \right) =\ker e_{1}^{\ast }.
\label{rowne2}
\end{equation}
Hence $\Psi \left( e_{1}^{\ast }\right) \circ S$ has a closed kernel, and
consequently is continuous for any $e_{1}^{\ast }\in E_{1}^{\ast }$. Since $%
\Psi $ is surjective, this means that $S$ becomes continuous when composed
with any continuous functional. So, by the Closed Graph Theorem $%
S:E_{1}\rightarrow E_{2}$ is continuous itself. We now can define a map $%
\widetilde{T}:\mathcal{A}_{2}\rightarrow B\left( E_{1}\right) $ by 
\begin{equation*}
\widetilde{T}\left( A\right) =S^{-1}\circ A\circ S\text{, for }A\in \mathcal{%
A}_{2}.
\end{equation*}
By (\ref{2}) 
\begin{equation*}
\widetilde{T}\circ T\left( e\otimes e^{\ast }\right) =S^{-1}\circ \left(
S\left( e\right) \otimes \Psi \left( e^{\ast }\right) \right) \circ
S=e\otimes \Psi \left( e^{\ast }\right) \circ S=e\otimes S^{\ast }\circ \Psi
\left( e^{\ast }\right)
\end{equation*}
and by (\ref{rowne2}) 
\begin{equation*}
\ker \Psi \left( e^{\ast }\right) \circ S=\ker e^{\ast },\text{ for all }%
e^{\ast }\in E_{1}^{\ast },
\end{equation*}
so each $e^{\ast }\in E_{1}^{\ast }$ is an eigenvector of $S^{\ast }\circ
\Psi ,$ hence $S^{\ast }\circ \Psi =\alpha Id$, and 
\begin{equation}
\widetilde{T}\circ T\left( e\otimes e^{\ast }\right) =\alpha e\otimes
e^{\ast }\text{,\quad for all }e\in E_{1},e_{1}^{\ast }\in E^{\ast }.
\label{3}
\end{equation}
Since $\widetilde{T}\circ T$ is linear it follows that 
\begin{equation*}
\widetilde{T}\circ T\left( A\right) =\alpha A
\end{equation*}
for all finite dimensional operators $A$. Notice that $\alpha \neq 0$ since
both $T$ and $\widetilde{T}$ are invertible.

Assume that $\widetilde{T}\circ T\neq \alpha Id$ and let $A_{0}\in \mathcal{A%
}_{1},$ $\widetilde{e}\in E_{1}$ be such that $\widetilde{T}\circ T\left(
A_{0}\right) \neq \alpha A_{0}$ and $\alpha A_{0}\left( \widetilde{e}\right)
\neq \left( \widetilde{T}\circ T\left( A_{0}\right) \right) \left( 
\widetilde{e}\right) .$ Put 
\begin{equation*}
B_{0}=A_{0}-A_{0}\left( \widetilde{e}\right) \otimes \widetilde{e}^{\ast },
\end{equation*}
where $\widetilde{e}^{\ast }\in E_{1}^{\ast }$ is such that $\widetilde{e}%
^{\ast }\left( \widetilde{e}\right) =1$. We have 
\begin{equation*}
B_{0}\circ \left( \widetilde{e}\otimes \widetilde{e}^{\ast }\right) =\left(
A_{0}-A_{0}\left( \widetilde{e}\right) \otimes \widetilde{e}^{\ast }\right)
\left( \widetilde{e}\right) \otimes \widetilde{e}^{\ast }=0,
\end{equation*}
while by (\ref{3}) 
\begin{align*}
\left( \widetilde{T}\circ T\left( B_{0}\right) \right) \circ \left( 
\widetilde{T}\circ T\left( \widetilde{e}\otimes \widetilde{e}^{\ast }\right)
\right) & =\alpha \left( \widetilde{T}\circ T\left( B_{0}\right) \right)
\circ \left( \widetilde{e}\otimes \widetilde{e}^{\ast }\right) \\
& =\alpha \left( \widetilde{T}\circ T\left( B_{0}\right) \right) \left( 
\widetilde{e}\right) \otimes \widetilde{e}^{\ast } \\
& =\alpha \left( \widetilde{T}\circ T\left( A_{0}\right) -\alpha \left(
A_{0}\right) \left( \widetilde{e}\right) \otimes \widetilde{e}^{\ast
}\right) \left( \widetilde{e}\right) \otimes \widetilde{e}^{\ast } \\
& =\alpha \left( \widetilde{T}\circ T\left( A_{0}\right) \left( \widetilde{e}%
\right) -\alpha \left( A_{0}\right) \left( \widetilde{e}\right) \right)
\otimes \widetilde{e}^{\ast } \\
& \neq 0
\end{align*}
which is a contradiction since $\widetilde{T}\circ T$ is biseparating. Hence 
$\widetilde{T}\circ T=\alpha Id$ so 
\begin{equation*}
TA=\alpha \widetilde{T}^{-1}A=\alpha S\circ A\circ S^{-1}\text{, \quad for }%
A\in \mathcal{A}_{1}.
\end{equation*}
\end{proof}

The next theorem extends the last result to the algebras of continuous
operator valued functions. For a Hausdorff set $X$ and a normed algebra $%
\mathcal{A}$ we denote by $C\left( X,\mathcal{A}\right) $ the algebra of all
continuous $\mathcal{A}$-valued functions on $X$ with the obvious
multiplication defined by 
\begin{equation*}
F\cdot G\left( x\right) =F\left( x\right) \cdot G\left( x\right) \text{%
,\quad\ for }x\in X\text{ and }F,G\in C\left( X,\mathcal{A}\right) ;
\end{equation*}
$C_{b}\left( X,\mathcal{A}\right) $ is the subalgebra of $C\left( X,\mathcal{%
A}\right) $ consisting of norm bounded functions. We equip $C_{b}\left( X,%
\mathcal{A}\right) $ with the sup norm topology and $C\left( X,\mathcal{A}%
\right) $ with the topology of uniform convergence on compact subsets of $X$%
. We denote by $C\left( X\right) $ the Banach algebra of all scalar valued
continuous functions on $X$ and by $C^{-1}\left( X\right) $ the set of
invertible elements of $C\left( X\right) $; we use notation $C_{b}\left(
X\right) $ for the Banach algebra of all bounded scalar valued continuous
functions on $X.$

For a function $f\in C\left( X\right) $ it may be often convenient to extend 
$f$ to a continuous function on $\beta X$, the maximal compactification of $%
X.$ In general, for $x\in\beta X\backslash X,$ the value of $f\left(
x\right) $ may be infinite. There are however completely regular spaces $X$
with a point $x_{0}\in\beta X\backslash X$ such that the value of $f\left(
x_{0}\right) $ is finite for all $f\in C\left( X\right) $; the set of all
points in $\beta X$ with this property is called the realcompactification of 
$X$ and is denoted by $rX.$ Since we have $C\left( X\right) =C\left(
rX\right) $ the natural domain for a continuous function on $X$ is $rX$, not 
$X$. Hence we will often consider realcompactifications, or alternatively we
will assume that the completely regular spaces under consideration are
realcompact. All compact sets are clearly realcompact, also all Lindel\"{o}f
spaces, and all metrizable spaces of nonmeasurable cardinal are realcompact
( \cite{GJ}, p.232).

For Banach spaces $E_{1},E_{2},$ $B\left( E_{1},E_{2}\right) $ is the space
of all linear continuous maps from $E_{1}$ into $E_{2}$ equipped with the
norm topology and $B^{-1}\left( E_{1},E_{2}\right) $ is the subset (possibly
empty) of $B\left( E_{1},E_{2}\right) $ consisting of invertible
isomorphisms.

Since $C\left( X,\mathcal{A}\right) $ is an algebra, Definition 1 describes
the meaning of a separating map from $C\left( X_{1},\mathcal{A}_{1}\right) $
into $C\left( X_{2},\mathcal{A}_{2}\right) $. However for the spaces of
vector valued continuous functions there is a possible alternative but not
equivalent natural definition: here we will call $T:C\left( X_{1},\mathcal{A}%
_{1}\right) \rightarrow C\left( X_{2},\mathcal{A}_{2}\right) $ \emph{%
strictly separating} if 
\begin{equation*}
\left\| F_{1}\left( \cdot\right) \right\| \left\| F_{2}\left( \cdot\right)
\right\| \equiv0\Longrightarrow\left\| \left( T\left( F_{1}\right) \right)
\left( \cdot\right) \right\| \left\| \left( T\left( F_{2}\right) \right)
\left( \cdot\right) \right\| \equiv 0\text{,\quad\ for }F_{1},F_{2}\in
C\left( X_{1},\mathcal{A}_{1}\right) .
\end{equation*}
That definition can be applied also if $\mathcal{A}$ is not an algebra but
just a normed linear space. We will later refer to \cite{AJ2} where a
general form of a strictly biseparating map is given.

\begin{theorem}
\label{th2}Assume $E_{1},E_{2}$ are Banach spaces and $\mathcal{A}_{1},%
\mathcal{A}_{2}$ are standard subalgebras of $B\left( E_{1}\right) $ and of $%
B\left( E_{2}\right) $, respectively. Assume further that $X_{1},X_{2}$ are
realcompact spaces and $T:C\left( X_{1},\mathcal{A}_{1}\right) \rightarrow
C\left( X_{2},\mathcal{A}_{2}\right) ,$ or $T:C_{b}\left( X_{1},\mathcal{A}%
_{1}\right) \rightarrow C_{b}\left( X_{2},\mathcal{A}_{2}\right) $ is a
biseparating map. Then $T$ is continuous and of the form 
\begin{equation*}
\left( T\left( F\right) \right) \left( x\right) =\alpha \left( x\right)
S_{x}\circ F\left( \varphi \left( x\right) \right) \circ S_{x}^{-1}\text{%
,\quad\ for }x\in X_{2}\text{ and }F\in \mathrm{domain}\left( T\right) ,
\end{equation*}
where $x\longmapsto S_{x}$ is a continuous map from $X_{2}$ into $%
B^{-1}\left( E_{1},E_{2}\right) ,$ $\varphi $ is a homeomorphism from $X_{2}$
onto $X_{1}$, and $\alpha \in C^{-1}\left( X_{2}\right) $.
\end{theorem}

\begin{proof}
Again we first need to introduce some notation. For $i=1,2$ we denote $%
\widetilde{C}\left( X_{i},\mathcal{A}_{i}\right) $ either $C\left( X_{i},%
\mathcal{A}_{i}\right) $ or $C_{b}\left( X_{i},\mathcal{A}_{i}\right) $
depending on the domain of the map $T;$ for $F\in \widetilde{C}\left( X_{i},%
\mathcal{A}_{i}\right) $ we put 
\begin{align*}
c\left( F\right) & =\{x\in X_{i}:F\left( x\right) \neq 0\}, \\
L\left( F\right) & =\{G\in \widetilde{C}\left( X_{i},\mathcal{A}_{i}\right)
:G\cdot F=0\}, \\
R\left( F\right) & =\{G\in \widetilde{C}\left( X_{i},\mathcal{A}_{i}\right)
:F\cdot G=0\}, \\
\mathcal{AI}_{i}& =\{H\in \widetilde{C}\left( X_{i},\mathcal{A}_{i}\right)
:L\left( H\right) \subset R\left( H\right) \}, \\
\mathcal{C}\left( F\right) & =\{G\in \widetilde{C}\left( X_{i},\mathcal{A}%
_{i}\right) :\forall H\in \mathcal{AI}_{i}\text{ }[F\cdot H=0\Longrightarrow
G\cdot H=0]\}.
\end{align*}
Notice that the sets $L\left( F\right) ,R\left( F\right) ,\mathcal{AI}_{i}$,
and $\mathcal{C}\left( F\right) $ have been defined solely using the
properties that are preserved by $T$ hence 
\begin{equation*}
T\left( L\left( F\right) \right) =L\left( T\left( F\right) \right) ,T\left(
R\left( F\right) \right) =R\left( T\left( F\right) \right) ,T\left( \mathcal{%
AI}_{1}\right) =\mathcal{AI}_{2},T\left( \mathcal{C}\left( F\right) \right) =%
\mathcal{C}\left( T\left( F\right) \right) .
\end{equation*}

We show that 
\begin{equation}
F_{1}\cdot F_{2}=0\Longleftrightarrow c\left( F_{1}\right) \cap c\left(
F_{2}\right) =\emptyset \text{,\quad\ for }F_{1}\in \widetilde{C}\left(
X_{i},\mathcal{A}_{i}\right) ,F_{2}\in \mathcal{AI}_{i}.  \label{sep}
\end{equation}
The implication $\Longleftarrow $ is obviously true for all functions in $%
\widetilde{C}\left( X_{i},\mathcal{A}_{i}\right) $. Assume $F_{1}\cdot
F_{2}=0$ and $x_{0}\in c\left( F_{1}\right) \cap c\left( F_{2}\right) $.
Since both $F_{1}\left( x_{0}\right) $ and $F_{2}\left( x_{0}\right) $ are
nonzero maps there is a continuous one dimensional linear map $A$ on $E_{i}$
such that $A\left( F_{1}\left( x_{0}\right) \right) \notin \ker F_{2}\left(
x_{0}\right) .$Put 
\begin{equation*}
G\left( x\right) =A\circ F_{1}\left( x\right) \text{,\quad\ for }x\in X_{i}%
\text{.}
\end{equation*}
We have 
\begin{equation*}
G\cdot F_{2}=0\text{ while }F_{2}\cdot G\left( x_{0}\right) =F_{2}\left(
x_{0}\right) \circ \left( A\left( F_{1}\left( x_{0}\right) \right) \right)
\neq 0
\end{equation*}
so $F_{2}\notin \mathcal{AI}_{i},$ which concludes the proof of (\ref{sep}).

By (\ref{sep}) 
\begin{align*}
\mathcal{C}\left( F\right) & =\{G\in \widetilde{C}\left( X_{i},\mathcal{A}%
_{i}\right) :\forall H\in \mathcal{AI}_{i}\text{ }[c\left( F\right) \cap
c\left( H\right) =\emptyset \Longrightarrow c\left( G\right) \cap c\left(
H\right) =\emptyset ]\} \\
& =\{G\in \widetilde{C}\left( X_{i},\mathcal{A}_{i}\right) :c\left( G\right)
\subset int\left( \overline{c\left( F\right) }\right) \},
\end{align*}
for arbitrary open sets $K,L$ we have $K\cap L=\emptyset $ if and only if $%
int\left( \overline{K}\right) \cap int\left( \overline{L}\right) =\emptyset $%
, so we get 
\begin{equation*}
c\left( F_{1}\right) \cap c\left( F_{2}\right) =\emptyset
\Longleftrightarrow \mathcal{C}\left( F_{1}\right) \cap \mathcal{C}\left(
F_{2}\right) =\{0\}\text{,\quad\ for }F_{1},F_{2}\in \widetilde{C}\left(
X_{i},\mathcal{A}_{i}\right) .
\end{equation*}
Since $T\left( \mathcal{C}\left( F\right) \right) =\mathcal{C}\left( T\left(
F\right) \right) $, for $F\in C\left( X_{i},\mathcal{A}_{i}\right) $, the
above proves that $T$ is strictly biseparating: 
\begin{equation}
\left\| F_{1}\left( \cdot \right) \right\| \left\| F_{2}\left( \cdot \right)
\right\| \equiv 0\Longleftrightarrow \left\| \left( T\left( F_{1}\right)
\right) \left( \cdot \right) \right\| \left\| \left( T\left( F_{2}\right)
\right) \left( \cdot \right) \right\| \equiv 0\text{,\quad\ for }%
F_{1},F_{2}\in \widetilde{C}\left( X_{i},\mathcal{A}_{i}\right) .  \label{ss}
\end{equation}
We need the following result from \cite{AJ2}.

\begin{theorem}
\label{ja}Assume $N_{1},N_{2}$ are normed spaces, $X_{1},X_{2}$ are
realcompact spaces, and $T:\widetilde{C}\left( X_{1},N_{1}\right)
\rightarrow \widetilde{C}\left( X_{2},N_{2}\right) $ is a linear bijection
satisfying (\ref{ss}). Then there is a bijective homeomorphism $\varphi
:X_{2}\rightarrow X_{1}$ and a map $J$ from $X_{2}$ into the set of linear
bijection from $N_{1}$ onto $N_{2}$ such that 
\begin{equation}
\left( T\left( F\right) \right) \left( x\right) =\left( J\left( x\right)
\right) \left( F\left( \varphi \left( x\right) \right) \right) ,\quad \ 
\text{for }x\in X_{2},\text{ and }F\in \widetilde{C}\left(
X_{1},N_{1}\right) .  \label{jaw}
\end{equation}
\end{theorem}

To finish the proof of Theorem \ref{th2}\ let $A_{1},A_{2}\in \mathcal{A}%
_{1} $ be such that $A_{1}\circ A_{2}=0,$ and denote by $\mathbf{A}_{1},%
\mathbf{A}_{2}$ the constant functions on $X_{2}$ equal to $A_{1}$, and to $%
A_{2}$, respectively. By Theorem \ref{ja} for any $x\in X_{2}$ 
\begin{equation*}
0=T\left( \mathbf{A}_{1}\right) \cdot T\left( \mathbf{A}_{2}\right) \left(
x\right) =\left( J\left( x\right) \right) \left( A_{1}\right) \circ \left(
J\left( x\right) \right) \left( A_{2}\right) ,
\end{equation*}
so $J\left( x\right) :\mathcal{A}_{1}\rightarrow \mathcal{A}_{2}$ is
separating. By the same arguments applied to $T^{-1}$ we conclude that $%
J\left( x\right) $ is biseparating. By Theorem \ref{th1} 
\begin{equation*}
\left( J\left( x\right) \right) \left( A\right) =\alpha _{x}S_{x}\circ
A\circ S_{x}^{-1}\text{,\quad\ for }A\in \mathcal{A}_{1},x\in X_{2},
\end{equation*}
where $S_{x}\in B^{-1}\left( E_{1},E_{2}\right) $. Hence, by (\ref{jaw}) 
\begin{equation*}
\left( T\left( F\right) \right) \left( x\right) =\alpha _{x}S_{x}\circ
\left( F\left( \varphi \left( x\right) \right) \right) \circ
S_{x}^{-1},\quad \ \text{for }x\in X_{2},\text{ and }F\in \widetilde{C}%
\left( X_{1},\mathcal{A}_{1}\right) ,
\end{equation*}
to check that $x\longmapsto \alpha _{x}$ is a continuous function it is
enough to put into the above formula $F$ equal, at every point of $X_{1},$
to the identity map on $E_{1}$.
\end{proof}

\begin{conjecture}
Assume $T$ is a biseparating map between $C^{\ast}$ algebras $\mathcal{A}%
_{1} $ and $\mathcal{A}_{2}$. Then $T$ is continuous and 
\begin{equation*}
T\left( a\right) =A_{0}\circ\Psi\left( a\right) \text{, \quad for }a\in%
\mathcal{A}_{1},
\end{equation*}
where $A_{0}$ is in the commutant of $\mathcal{A}_{2}$.
\end{conjecture}

\end{document}